\documentclass[10pt]{article}

\usepackage[top=4cm, bottom=3cm, left=3cm, right=3cm]{geometry}
\usepackage{graphicx}
\usepackage[parfill]{parskip}
\usepackage{amsmath, amsthm, amssymb,amsfonts}
\usepackage{wasysym}
\usepackage[utf8]{inputenc}
\usepackage[english, frenchb]{babel}
\usepackage[T1]{fontenc}
\usepackage{url}
\usepackage{verbatim}
\usepackage{indentfirst}
\usepackage{xypic}
\usepackage{braket}
\usepackage{fancyhdr}
\usepackage{enumerate}
\usepackage{textcomp}
\usepackage{stmaryrd}
\usepackage{xcolor}
\usepackage[colorlinks=true, linktoc=page, citecolor=blue, linkcolor=blue, urlcolor=blue]{hyperref}
\usepackage{fancyhdr}

\usepackage{mathptmx}
\usepackage{mathrsfs}
\usepackage{lipsum}
\usepackage{colortbl}

\DeclareFontFamily{OT1}{pzc}{}
\DeclareFontShape{OT1}{pzc}{m}{it}{<-> s * [1.1] pzcmi7t}{}
\DeclareMathAlphabet{\mathpzc}{OT1}{pzc}{m}{it}  

\usepackage{setspace}
\usepackage{titlesec}
\renewcommand\thepart{\textbf{Part \Roman{part}.}}
\titleformat{\part}
  {\sc  \center}
  {\thepart}{1em}{}

\titleformat{\section}
  {\bf \normalsize \center}
  {\thesection.}{1em}{}

\titleformat{\subsection}
  {\bf \normalsize}
  {\thesubsection.}{1em}{}

\titleformat{\subsubsection}
  {\bf \normalsize}
  {\thesubsubsection.}{1em}{}  %
\newtheoremstyle{defn}{}{}{}{}{\bfseries}{ : }{0cm}{}
\newtheoremstyle{theoreme}{6pt}{6pt}{\itshape}{}{\bfseries}{.}{3pt}{}

\theoremstyle{remark} 
\newtheorem*{exemple}{Example} 

\theoremstyle{defn} 

\theoremstyle{theoreme} 
\newtheorem{thm}{Theorem}[section]

\newtheorem{lemme}[thm]{Lemma}

\newtheorem{prop}[thm]{Proposition}

\newtheorem{cor}[thm]{Corollary}
\newtheorem*{def-prop}{Définition-Proposition}

\makeatletter
\renewenvironment{proof}[1][Proof:]{\par
\normalfont
\topsep0\p@\@plus6\p@ \trivlist
\item[\hskip\labelsep\itshape
#1]
}{%
\qed\\ \endtrivlist
}
\makeatother

\numberwithin{equation}{section}

\newcounter{num} 
\newcommand\exo{\stepcounter{num}\noindent{\bf Exercice~\thenum} \kern10pt }

\newcommand{\R}{\mathbb{R}}
\newcommand{\C}{\mathbb{C}}

\newcommand{\N}{\mathbb{N}}

\newcommand{\g}{\mathfrak{g}}

\newcommand{\norme}[1]{\left\Vert #1\right\Vert}

\newcommand{\ind}{\text{Ind}}

\title{\hspace{-0.2cm}On~the~analogy~between~real~reductive~groups and Cartan motion groups. III:~A~proof~of~the~Connes-Kasparov~isomorphism}
\author{Alexandre Afgoustidis~\\ \small \emph{CEREMADE, Université Paris-Dauphine}}
\date{}

\makeatletter
\def\@settitle{\begin{center}%
  \baselineskip14\p@\relax
    \bfseries
    \normalfont\huge
  \@title
  \end{center}%
}
\def\@settitle{\begin{center}%
  \baselineskip14\p@\relax
    \bfseries
    \normalfont\huge
  \@title
  \end{center}%
}
\makeatother

\begin{document}
\selectlanguage{english}
\maketitle
\begin{abstract} \noindent Alain Connes and Nigel Higson pointed out in the 1990s that the Connes-Kasparov ``conjecture'' for the K-theory of reduced group $C^\ast$-algebras seemed, in the case of reductive Lie groups, to be a cohomological echo of a conjecture of George Mackey concerning the rigidity of representation theory along the deformation from a real reductive group to its Cartan motion group. For complex semisimple groups, Nigel Higson established in 2008 that Mackey's analogy is a real phenomenon, and does lead to a simple proof of the Connes-Kasparov isomorphism. We here turn to more general reductive groups and use our recent work on Mackey's proposal, together with Higson's work, to obtain a new proof of the Connes-Kasparov isomorphism.  \end{abstract}

\section{Introduction}\label{subsec:conneskasp_intro} When $G$ is a second countable locally compact group, the Fell topology on its unitary dual $\widehat{G}$ is in general quite wild: studying the topological space $\widehat{G}$ directly is usually difficult. An indirect approach, that has often proved fruitful, is to study a suitable completion of the convolution algebra $C^\infty_{c}(G)$ of continuous and compactly supported functions on $G$. The various completions will be noncommutative algebras, and may be thought of as noncommutative-geometry replacements for the (not very helpful) space of continuous functions on $\widehat{G}$.

We shall be concerned with the \emph{reduced $C^\star$-algebra} of $G$. Assume the group $G$ to be unimodular and equipped with a Haar measure. When $f$ is an element of $C^\infty_{c}(G)$, convolution with $f$ defines a bounded linear operator on the Hilbert space $\mathbf{L}^2(G)$, and that operator has a norm, say $\norme{f}$; the {reduced} $C^\star$-algebra $C^\star_{r}(G)$ is the completion of $C^\infty_{c}(G)$ with respect to $\norme{\cdot}$. The spectrum of $C^\star_{r}(G)$ is the \emph{reduced dual} $\widetilde{G}$; it gathers the unitary irreducible representations that are, loosely speaking, necessary to ``decompose'' the regular representation of $G$ on $\mathbf{L}^2(G)$ into irreducibles.

The \emph{Baum-Connes conjecture} describes the \emph{$K$-theory} of $C^\star_{r}(G)$, to be thought of as a ``non-commutative'' replacement for the Atiyah-Hirzebruch K-theory of $\widetilde{G}$ (see \cite{BaumConnesHigson} for the formulation of the conjecture, and \cite{KasparovICM} for the K-theory). Although we will consider only reductive groups and their Cartan motion groups in this note, we should mention that the interest of the conjecture partly stems from its generality (it encodes very nontrivial features of both Lie groups and discrete groups, and the existence of an analogue for groupoids makes it suitable for the study of foliations), and its deep connections with geometry, index theory and topology.

\paragraph{The Connes-Kasparov isomorphism for connected Lie groups.}  When $G$ is a connected \emph{Lie group}, the conjecture is equivalent with the assertion that the reduced dual of $G$ can be accounted for, at least at the level of K-theory, with the help of Dirac operators. Suppose $K$ is a maximal compact subgroup of $G$, and $R(K)$ is the representation ring of $K$ $-$ whose underlying abelian group is freely generated by the equivalence classes of irreducible $K$-modules. Starting from an irreducible $K$-module, and after going up to a two-fold covering of $G$ and $K$ if necessary, one can build an equivariant spinor bundle on $G/K$ and a natural $G$-invariant elliptic operator (the Dirac operator) acting on its sections; that operator has an index which can be refined into an element of  the K-theory group $\mathpzc{K}_j\left[C^\star(G)\right]$, where $j \equiv \dim(G/K)$ mod $2$. This produces a morphism of abelian groups
\begin{equation} \label{ConnesKasparov} R(K) \overset{\mu}{\longrightarrow} \mathpzc{K}_j\left[C^\star_r(G)\right] \quad (j \equiv  \dim(G/K) \ \text{mod} \ 2)\end{equation}
called Dirac induction; the \emph{Connes-Kasparov conjecture} is the statement that $\mu$ is an isomorphism, and that the K-theory group $\mathpzc{K}_j\left[C^\star_r(G)\right]$ is zero when $j \not\equiv \dim(G/K)$ mod $2$. See \cite[\S 4]{BaumConnesHigson}.

When $G$ is reductive, our current understanding of the reduced (or: tempered) dual $\widetilde{G}$ is in many respects complete. As a topological space, $\widetilde{G}$ is in addition much more reasonable than the reduced duals studied by the full-blown Baum-Connes conjecture: it is (roughly speaking) akin to a real affine variety. Yet the reductive case is an important one: a major source for the formulation of the Baum-Connes conjecture was the Partasarathy-Atiyah-Schmid realization of the discrete series using Dirac operators (see \cite{Parthasarathy, AtiyahSchmid}), and understanding the reductive case proved to be the key to understanding the case of general Lie groups \cite{ChabertEchterhoffNest}.

\paragraph{Previous proofs for reductive Lie groups.} When $G$ is a linear connected reductive group, the Connes-Kasparov ``conjecture'' has been a classical result for thirty years. Two proofs have been given, and they are quite different:

$\bullet$ Antony Wassermann's short note of 1987 \cite{Wassermann} uses the comprehensive knowledge of $C^\star_{r}(G)$ extracted by Arthur \cite{ArthurSchwartz} from Harish-Chandra's  work; his proof consists in an explicit calculation of the right-hand side of \eqref{ConnesKasparov} and the arrow therein. This followed earlier treatment of special cases from the same perspective: the important but simpler case of complex semisimple groups had been covered by Penington and Plymen \cite{PeningtonPlymen}, and Valette had considered two classes of real groups (see \cite{ValetteRang1, Valette2}). Wassermann did not provide a more detailed exposition of his argument, which uses the Knapp-Stein theory of intertwining operators. Clare, Crisp and Higson recently gave a quite accessible account of the structure of $C^\star_r(G)$  along Arthur's lines \cite{ClareCrispHigson}, and should soon provide related insight on the K-theoretic side.

$\bullet$ Vincent Lafforgue developed in 1998 a deep notion of bivariant  KK-theory for group actions on Banach spaces (\cite{Lafforgue}; see also \cite{SkandalisBourbakiLafforgue}). This opened him a way to the Baum-Connes isomorphism that is not only very well-suited to reductive Lie groups, but also encompasses reductive $p$-adic groups as well as some discrete subgroups which had resisted every approach before his. His strategy is almost orthogonal to that of Wassermann,  replacing most of the arsenal of representation theory by a few simple (but far-reaching) facts on the distance to the origin in $G/K$ and on Harish Chandra's elementary spherical functions. His framework proved flexible enough to lend itself to the extensions needed to prove the Connes-Kasparov conjecture when $G$ is an \emph{arbitrary} Lie group (with a finite number of connected components): Chabert, Etcherhoff and Nest proved the Connes-Kasparov conjecture in 2003 \cite{ChabertEchterhoffNest} by using the decomposition of any such Lie group as a semidirect product of a reductive and a nilpotent group, and Mackey's theory for the representations of group extensions \cite{MackeyActa}.

\paragraph{Lie group deformations and the Mackey analogy.} A third way to the Connes-Kasparov conjecture has long been suspected to exist, should one dig deeply enough in the theory of group \emph{deformations} (or \emph{contractions} in the classical terminology of {I}n\"{o}n\"{u} and Wigner \cite{InonuWignerContractions}).

 Suppose $G$ is a connected  reductive Lie group, $K$ is a maximal compact subgroup, and write $G_0$ for the semidirect product $K \ltimes (\g/\mathfrak{k})$ associated to the adjoint action of $K$ on $\g/\mathfrak{k}$ (here $\g$ and $\mathfrak{k}$ denote the Lie algebras of $G$ and $K$). Then there is a "continuous" family of groups $(G_t)_{t \in \R}$ which interpolates between $G_1 := G$ and $G_0$ (see \S \ref{subsec:defo} below). This deformation gives rise to a continuous field $\left\{C^\star(G_t)\right\}_{t \in \R}$ of algebras; Alain Connes and Nigel Higson  observed in the late 1980s and early 1990s (see \cite{ConnesHigsonDeformations}, \cite[\S 4]{BaumConnesHigson}, \cite[\S II.10.$\beta$]{ConnesLivreENG}) that the Connes-Kasparov conjecture is equivalent with the fact that this field has \emph{constant} K-theory.
 
 As Connes and Higson insisted, this meant that the Connes-Kasparov isomorphism could be the noncommutative-geometric counterpart of an intriguing rigidity phenomenon occurring at the level of representation theory; in fact, the reformulation of the Connes-Kasparov conjecture in terms of deformations echoed enthusiastic observations by G.~W. Mackey \cite{MackeyConjecture} on a possible relationship between the representation theories of $G$ and $G_0$ when $G$ is a semisimple Lie group. Mackey had conjectured that  in the reductive case, there were deep-rooted, though surprising, analogies between $\widetilde{G}$ and $\widehat{G_0}=\widetilde{G_0}$. Connes and Higson's observations strongly invited to view the isomorphism \eqref{ConnesKasparov} a simple K-theoretic reflection of these analogies, roughly indicating that $\widetilde{G}$ and $\widetilde{G_0}$ share algebraic-topological invariants (see \cite{BaumConnesHigson}, \S 4). As Connes put it  (\cite{ConnesLivreENG}, $\S$ II.10.$\delta$):\\

\begin{center}
\begin{minipage}{0.95\textwidth}
\emph{It is of course desirable to find direct proofs of the surjectivity of [the Baum-Connes assembly map]. In that respect the ideas developed by Mackey in {\cite{MackeyConjecture}}, or in the theoretical physics literature on deformation theory, should be relevant.\\ }
\end{minipage}
\end{center}

 
\noindent Connes and Higson had crafted precise tools to do so in 1990, to be recalled below (see \S \ref{subsec:defo}). But Mackey's conjecture lay dormant for a long time; only recently did it become clear that the strategy imagined by Connes and Higson can be pursued for any reductive group. Going through that path is the aim of the present~note.


\paragraph{Higson's work on complex semisimple groups.} About ten years ago, Higson decided to re-examine the connections between Mackey's observations and K-theory in the case of complex semisimple groups \cite{Higson2008}. For that special class of reductive groups, he proved the existence of a natural bijection between the reduced duals $\widetilde{G}$ and $\widetilde{G_0}$; the representation-theoretic properties of the bijection revealed the structures of $C^\star_{r}(G_0)$ and $C^\star_{r}(G)$ to be close enough to lead to a surprisingly simple proof of  the Connes-Kasparov conjecture. When $G$ is a complex semisimple group, the following statement summarizes his analysis of the structure of the field $\left\{C^\star(G_t)\right\}_{t \in \R}$.
\vspace{-0.2cm}

\noindent \textbf{Theorem A.} \emph{ Consider the family $(G_t)_{t \in \R}$ of Lie groups. \begin{enumerate}[(i)]
\item For every $t$ in $\R$ (including $t = 0$), the partition of $\widetilde{G_t}$ according to lowest $K$-types determines a dense filtration $\mathbf{I}_t[1] \subset\mathbf{I}_t[2] \subset \hdots \subset\mathbf{I}_t[p] \subset \hdots \subset C^\ast_r(G_t)$ of the reduced $C^\ast$-algebra $C^\ast_r(G_t)$ by closed ideals.
The filtration at $t$ is the specialization of a dense filtration $\mathcal{I}[1] \subset \hdots \subset\mathcal{I}[p] \subset \hdots \subset \mathcal{C}$ of the $C^\ast$-algebra of sections of the continuous field $\left\{C^\star(G_t)\right\}_{t \in \R}$.
\item  For every $p$ in $\N$, the subquotient field $\mathcal{C}[p]=\mathcal{I}[p+1]/\mathcal{I}[p]$ is strongly Morita-equivalent with a \emph{constant} field (of commutative algebras).  \end{enumerate} }

\noindent Using part (i) and the cohomological nature of $K$-theory, the Connes-Kasparov isomorphism then appears as an easy consequence of the rigidity statement in (ii).

This way to the Connes-Kasparov conjecture does takes one through the fine structure of representation theory: the filtration is defined in representation-theoretic terms, and the rigidity statement in (ii) reflects a very peculiar behavior of matrix elements (of principal series representations) along the deformation.  But rather than using representation theory for a direct calculation of \eqref{ConnesKasparov}, it expresses the Connes-Kasparov phenomenon as a K-theoretic reflection of the structural analogy between $C^\star(G)$ and $C^\star(G_0)$ induced by the correspondence between $\widetilde{G}$ and $\widetilde{G_0}$. An appealing feature of this approach is that only simple and quite general facts about K-theory are needed, and that no K-theory group need be written down explicitly.

\paragraph{Contents of this note.} When $G$ is a real reductive group, we recently described a natural bijection between $\widetilde{G}$ and $\widetilde{G_0}$ \cite{AAMackey}. We will here prove that the Connes-Kasparov isomorphism for real reductive Lie groups can indeed be deduced in a relatively elementary manner from the representation-theoretic properties of that bijection, using the $C^\ast$-algebraic ideas and methods due to  Higson and Connes. 

The representation theory of general reductive Lie groups is more complicated than that of complex semisimple groups: for instance, the existence of the discrete series (and with it the need for the bulk of Harish-Chandra's work) is a specific feature of the real case. The proof of Theorem A in Higson's work on complex groups accordingly rests on representation-theoretic facts which, on the surface, may look quite special to the complex case. It may therefore seem surprising that Theorem A above should be obtained for real reductive groups in much the same way as it can for complex groups; yet that is exactly what will happen. We will see that the difficulties that one could expect to occur for real groups, ranging from the non-uniqueness of lowest $K$-types to the necessity of tracing the collapse of some matrix elements of discrete series representations (of reductive subgroups) along the deformation, are reflected precisely enough in the properties of the Mackey-Higson bijection for real groups that the complications cancel out as the contraction is performed. Thus, the present work will not necessitate new ingredients on the $C^\star$-algebra side $-$  most ideas in \cite{Higson2008} (and a great many technical lemmas) will turn out to be usable without any conceptual change, and it should be very clear that all steps of our proof are either due to Higson or strongly inspired by his work. The properties of the Mackey-Higson bijection recalled in \S \ref{subsec:rappels_mackey} do of course draw from several important sources, which include the Knapp-Zuckerman classification of irreducible tempered representations and Vogan's work on lowest $K$-types.

\paragraph{Acknowledgments.} I am very grateful to Daniel Bennequin, Nigel Higson and Georges Skandalis for their help and support. A first version of this work appeared in Chapter 8 of my Ph.D. thesis \cite{AAthese}, prepared at Université Paris-7 under the guidance of Daniel Bennequin. 

\section{Background and notations}

\subsection{The deformation field; Connes and Higson's deformation interpretation of the assembly map \eqref{ConnesKasparov}} \label{subsec:defo}

\noindent Let $G$ be a linear connected reductive group. We fix a maximal compact subgroup $K$, write $\g$ and $\mathfrak{k}$ for the Lie algebras of $G$ and $K$, and let $\g = \mathfrak{k} \oplus \mathfrak{p}$ be the corresponding Cartan decomposition of $\g$. The \emph{Cartan motion group} $G_0$ is the semidirect product $G_0 = K \ltimes \mathfrak{p}$ associated with the adjoint~action~of~$K$~on~$\mathfrak{p}$. 

For every nonzero real number $t$, we define a group $G_t$ by using the global diffeomorphism 
\begin{align*} \varphi_t :  K \times \mathfrak{p} & \to G \\ (k,v)&  \mapsto k\exp_G(tv)\end{align*}
to endow the set $K \times \mathfrak{p}$ with the product law which turns $\varphi_t$ into an isomorphism. Thus, for every $t \in \R$, the group $G_t$ is \emph{equal}  as a topological space with $K \times \mathfrak{p}$. 

In order to give a precise meaning to the field of $C^\ast$-algebras to be considered below, we equip the disjoint~union \[ \mathpzc{G} := \bigsqcup \limits_{t \in \R} G_{t} \]
with the smooth manifold structure for which the bijection
\begin{align} \label{variables} \mathpzc{G} & \rightarrow G_{0} \cup \left( G \times \R^\times \right)\\ \nonumber \gamma \in G_t  & \mapsto \begin{cases} \gamma & \text{ \ if $t=0$} \\ (\varphi_{t}(\gamma),t) & \text{ \ if $t \neq 0$ } \end{cases}\end{align}
is turned into a diffeomorphism when $G_{0} \cup \left( G \times \R^\times \right)$ is equipped with the smooth structure described in section 6.2 of \cite{Higson2008}. There is a natural smooth family of measures on $K \times \mathfrak{p}$ giving a Haar measure for each group $G_t$, $t \in \R$.  We can consider the reduced $C^\star$ algebra of each $G_{t}$ for the corresponding Haar measure;  the field
\[ \left\{ {C}_{r}^\star(G_{t}) \right\}_{t \in \R} \]
is then a continuous field of $C^\star$-algebras (with our choice of smooth structure on $\mathpzc{G}$, this follows from Lemma 6.13 in \cite{Higson2008}). We will consider 
\begin{equation} \label{champ_c} \mathcal{C}:= \text{ $C^\star$ algebra of continuous sections of the restriction of $\left\{ {C}_{r}^\star(G_{t}) \right\}$ to  [0,1].  }\end{equation}
 The evaluation maps at $t=0$ and $t=1$ induce $C^\star$-algebra morphisms from $\mathcal{C}$ to $C^\star_{r}(G_0)$ and $C^\star_{r}(G)$, respectively, and in turn these induce two homomorphisms $\alpha_{0}: \mathpzc{K}(\mathcal{C}) \rightarrow \mathpzc{K}(C_{r}^\star(G_0))$ and $\alpha_{1}: \mathpzc{K}(\mathcal{C}) \rightarrow \mathpzc{K}(C_{r}^\star(G_1))$. Because the field $\left\{ C^\star(G_t) \right\}_{t \in ]0,1]}$ (with $t=0$ excluded) is trivial, $\alpha_0$ is an isomorphism.  Now, the composition
\begin{equation} \label{CH} \alpha_1 \circ \alpha_0^{-1}: \mathpzc{K}(C^\star(G_0)) \rightarrow \mathpzc{K}(C^\star(G))\end{equation}
has an important connection with the Connes-Kasparov conjecture (see \cite[p. 24]{BaumConnesHigson}). Because $G_0 = K \ltimes \mathfrak{p}$ is an extension of $K$ by a vector abelian group, the equivariant Bott periodicity theorem yields a natural isomorphism between $\mathpzc{K}_j(C^\star(G_0))$ and $R(K)$ for $j \equiv \dim(G/K) \text{ mod } 2$, while $\mathpzc{K}_j(C^\star(G_0))=0$ for $j \not\equiv \dim(G/K) \text{ mod } 2$. Viewed through the isomorphism between $\mathpzc{K}_j(C^\star(G_0))$ and $R(K)$, the map $\alpha_1 \circ \alpha_0^{-1}$ in \eqref{CH} is none other than the assembly map \eqref{ConnesKasparov}:

\begin{thm}[Connes and Higson, 1990 \cite{ConnesHigsonDeformations}] \label{strategie_CK} The Connes-Kasparov assembly map \eqref{ConnesKasparov} is an isomorphism if and only if the map $\alpha_1$, induced by evaluation at $t = 1$, is an isomorphism too.\end{thm}

\subsection{Construction and properties of the Mackey-Higson correspondence} \label{subsec:rappels_mackey}

\paragraph{Lowest $K$-types.} Let $G$ be a linear reductive group with all Cartan subgroups abelian (see \cite[\S 0.1]{Vogan81}; if $G$ is linear connected reductive, then all its Cartan subgroups are abelian). Fix  a maximal compact subgroup $K$ in $G$, a Cartan subalgebra $\mathfrak{t}$ of the Lie algebra $\mathfrak{k}$ and a system $\Delta^+$ of positive roots for the pair $(\mathfrak{k}_\C, \mathfrak{t}_\C)$. Form the corresponding half-sum $\rho_c$ of positive roots; it is an element of the vector space dual $\mathfrak{t}^\star$ of $\mathfrak{t}$. Any choice of symmetric nondegenerate bilinear form on $\g$ that is negative-definite on $\mathfrak{k}$ and positive-definite on $\mathfrak{p}$ determines a Euclidean norm on $\mathfrak{t}^\star$.  For every class $\lambda$ in $\widehat{K}$, we write $\norme{\lambda}_{\widehat{K}}$ for the distance in $\mathfrak{t}^\star$ between the $\Delta^+$-highest weight of $\lambda$ and $-2\rho_c$. The magnitude function $\norme{\cdot}_{\widehat{K}}$ determines a partial order on $\widehat{K}$, which depends neither on the choice of $T$ nor on that of $\Delta^+$. 

Every irreducible tempered representation $\pi$ of ${G}$ or $G_0$ comes with a (finite) set of  \emph{lowest $K$-types}: these are the elements of $\widehat{K}$ that occur in the decomposition of the restriction $\pi_{|K}$ into irreducibles and are minimal (among the $K$-types that occur in $\pi_{|K}$) for the above ordering.

\begin{prop} \label{multiplicite} 
\begin{enumerate}[(a)]
\item If $\pi$ is an irreducible tempered representation of $G$, then every lowest $K$-type of $\pi$ occurs with multiplicity one in the restriction  $\pi_{|K}$.
\item If $\pi_0$ is a unitary irreducible representation of $G_0$, then every lowest $K$-type of $\pi_0$ occurs with multiplicity one in the restriction  $(\pi_0)_{|K}$.
\end{enumerate} \end{prop} 
Part (a) is a fundamental result of Vogan (proved in \cite{Vogan81}; see also \cite{Vogan77}); for (b), see \cite{AAMackey}, Corollary 4.3.

\paragraph{The Mackey-Higson bijection (see \cite{AAMackey}).}
 
We return to our linear connected reductive group $G$. Suppose $\chi \in \mathfrak{p}^\star$ is a linear functional on $\mathfrak{p}$ and $\mu$ is an irreducible representation of the stabilizer $K_\chi$ of $\chi$ in $K$ (for the coadjoint action of $K$ on $\mathfrak{p}^\star$). 
\begin{itemize}
\item[$\bullet$] The pair $(\chi, \mu)$ can be used to produce a unitary irreducible representation of $G_0$, the induced representation 
\[ \mathbf{M}_0(\chi, \mu)= \ind_{K_\chi \ltimes \mathfrak{p}}^{G_0}(\mu \otimes e^{i\chi}).\]
\item[$\bullet$] The same pair can be used to define a tempered irreducible representation of $G$, as well.

Write $L_\chi$ for the centralizer of $\chi$ in $G$ (for the coadjoint action). It turns out that there exists a parabolic subgroup $P_\chi = L_\chi N_\chi$ of $G$ with Levi factor $L_\chi$. Write $L_\chi = M_\chi A_\chi N_\chi$ for the Langlands decomposition of $P_\chi$. The  group $M_\chi$ admits $K_\chi$ as a maximal compact subgroup; it is usually disconnected, but all its Cartan subgroups are abelian. The representation $\mu$ of $K_\chi$ determines a unique irreducible tempered representation of $M_\chi$: among the tempered irreducible representations of $M_\chi$ with real infinitesimal character\footnote{An irreducible tempered representation has real infinitesimal character if and only if it occurs as an irreducible summand in a representation of the form $\ind_{MAN}^G(\tau \otimes \mathbf{1})$, where $MAN$ is a cuspidal parabolic subgroup of $G$ and $\tau$ is a discrete series representation of $M$.}, there is one and only one which admits $\mu$ as a lowest $K_\chi$-type (and then $\mu$ is its only lowest $K_\chi$-type). We write $\mathbf{V}_{M_\chi}(\mu)$ for it and define a representation of $G$ as 
\[ \mathbf{M}(\chi, \mu) = \ind_{M_\chi A_\chi N_\chi}^G(\mathbf{V}_{M_\chi}(\mu) \otimes e^{i\chi} \otimes \mathbf{1}).\]
That representation is irreducible and tempered \cite[Theorem 3.3(a)]{AAMackey}.
\end{itemize} 
Mackey's early work on semidirect products shows that every representation in $\widehat{G_0}=\widetilde{G_0}$ is of the form $\mathbf{M}_0(\chi, \mu)$ for some $(\chi, \mu)$; it turns out that the class of $\mathbf{M}(\chi, \mu)$ in $\widetilde{G}$ depends only on the class of $\mathbf{M}_0(\chi, \mu)$ in $\widehat{G_0}=\widetilde{G_0}$, so the correspondence $\mathbf{M}_0(\chi, \mu) \leftrightarrow \mathbf{M}(\chi, \mu)$  induces a map  $\mathcal{M}: \widetilde{G_0} \to \widetilde{G}$.

\begin{thm}[\cite{AAMackey}, Theorem 3.3] \label{correspondance} The map $\mathcal{M}: \widetilde{G_0} \to \widetilde{G}$ is a bijection between the reduced duals of $G_0$ and $G$. \end{thm}

\begin{prop}[\cite{AAMackey}, Proposition 4.1] \label{preservation_ktypes} The correspondence $\mathcal{M}: \widetilde{G_0} \to \widetilde{G}$ preserves lowest K-types. \end{prop}

\subsection{Notations on matrix coefficients} \label{subsec:notations_coeffs}

%

We recall some notations and elementary results which we will need to take up from \cite{Higson2008}. Suppose $G$ is a connected unimodular Lie group, $K$ is a compact subgroup, $s$ is a smooth function on $K$ and $f$ is a smooth and compactly supported function on $G$. Choose a Haar measure on $G$ and define two convolutions between $s$~and~$f$, 
\begin{equation} \label{convol_K} s \underset{K}{\star} f = g \mapsto \frac{1}{\text{Vol}(K)} \int_{K} s(k) f(k^{-1}g) dk \ \ , \ \ f \underset{K}{\star} s = g \mapsto \frac{1}{\text{Vol}(K)}  \int_{K} f(kg) s(k^{-1})dk. \end{equation}
These are two smooth and compactly supported functions on $G$.

 Now suppose $K$ is any compact Lie group, $K_1$ is a closed subgroup, $(V, \tau)$ is an irreducible representation of $K$ with orthonormal basis $\{ v_\alpha \}_{\alpha=1,\dots,\dim(V)}$, and $W$ is a $K_1$-invariant \emph{irreducible} subspace of $V$. Write $e_{\alpha \beta}$ for the matrix element $k \mapsto \dim(V) \langle \tau(k) v_\beta, v_\alpha\rangle$ (this is a smooth function on $K$). When $v_{\alpha}$ and $v_{\beta}$ both lie in $W$, write $d_{\alpha \beta}$ for $k \mapsto \dim(W) \langle \tau(k) v_\beta, v_\alpha\rangle$; the restriction of $d_{\alpha \beta}$ to $K_1$ is a matrix element of $(W, \tau|_{|K_1})$. Then the Schur-Weyl orthogonality relations yield
\[ e_{\alpha \beta} \underset{K}{\star} e_{\beta \gamma} = e_{\alpha \gamma} ;\]
\begin{equation} \label{weyl} d_{\alpha \beta} \underset{K_1}{\star} e_{\beta \gamma} = e_{\alpha \gamma} =  d_{\alpha \beta} \underset{K_1}{\star} e_{\beta \gamma} \end{equation} 
(for the second equality, it is to be assumed that $v_{\alpha}$, $v_{\beta}$ and $v_{\gamma}$ lie in $W$). We also note that $e_{\alpha \beta}(k) = \overline{e_{\beta \alpha}(k^{-1})}$ for all $k$ (the bar denotes complex conjugation).

\section{Distinguished subquotients of the reduced $C^\star$-algebras} \label{sec:sous_quotients}

\subsection{Subquotients of the reductive group's algebra} \label{sec:cstar}

%
%
%
%

\subsubsection{Ideals associated with a set of lowest $K$-types} \label{subsec:sous_quotient_G}

\noindent We return to studying our linear connected reductive group $G$ and maximal compact subgroup $K$. Define~a~set 
\[ \mathpzc{Classes} \subset \left\{\  \text{finite subsets of } \widehat{K}\right\} \]
by declaring that a finite subset $\mathpzc{C}$ of $\widehat{K}$ lies in $\mathpzc{Classes}$ when there exists an irreducible tempered representation of $G$ whose set of lowest $K$-types is exactly $\mathpzc{C}$. Note that in this case, $\norme{\cdot }_{\widehat{K}}$ (from \S \ref{subsec:rappels_mackey}) takes the same value on all the elements of $\mathpzc{C}$. 

 We will associate a subquotient of the reduced $C^*$-algebra ${C}_{r}^\star(G)$ to every class in $\mathpzc{Classes}$. Later on it will be convenient that the family subquotients obtained in this way be associated to an increasing sequence of ideals in ${C}_{r}^\star(G)$, so let us choose first a linear ordering
\[ \mathpzc{Classes} = \left\{ \mathpzc{C}_{1}, \mathpzc{C}_{2}, \dots\right\} \]
in such a way that 
\begin{itemize}
\item[$\bullet$] if the value of $\norme{\cdot }_{\widehat{K}}$ on $\mathpzc{C}_{p}$ is (strictly) smaller than that on  $\mathpzc{C}_{q}$, then $p < q$,
\item[$\bullet$] if the values of $\norme{\cdot }_{\widehat{K}}$ on $\mathpzc{C}_{p}$ and $\mathpzc{C}_{q}$ agree but the number of elements in $\mathpzc{C}_{p}$ is (strictly) larger than that in $\mathpzc{C}_{q}$, then $p < q$.

\end{itemize}

For every class  $\lambda$ in $\widehat{K}$, we fix an irreducible $K$-module  $V_{\lambda}$ with equivalence class $\lambda$, and write $\langle \cdot, \cdot \rangle$ for a $K$-invariant inner product on it. Any choice of nonzero vector $v$ in an irreducible $K$-module with class $\lambda$ determines a matrix element \begin{align} \label{choix_v} \mathbf{p}^{v}_{\lambda} : K & \to \C \\ k & \mapsto \langle v_{{}}, \lambda(k^{-1}) v_{{}} \rangle. \nonumber \end{align}
Convolution with this matrix element, as in \eqref{convol_K}, yields an element in the multiplier algebra of $  {C}_{r}^\star(G)$.  We can then define a closed ideal in ${C}_{r}^\star(G)$ by setting
\[ \mathbf{J}[p]  = \bigcap \limits_{\lambda \in \mathpzc{C}_{p}}  {C}_{r}^\star(G)  \mathbf{p}^v_{\lambda} {C}_{r}^\star(G),  \]
and obtain an increasing family $( \mathbf{I}[p] )_{p \in \N}$ of ideals by setting $\mathbf{I}[p]=\mathbf{J}[1] + \dots + \mathbf{J}[p]$ for each $p$. We will study the subquotients
\[ \mathbf{C}[{p}] = \mathbf{I}[p+1]/\mathbf{I}[p].\] 
The algebras $\mathbf{J}[p]$ and $\mathbf{C}[p]$ depend on a choice of vector $v=v_\lambda$ for each $\lambda \in \widehat{K}$, which will be made explicit in  \S \ref{subsec:explicite}. But their spectra do not depend on that choice. They are locally closed subsets of $\widehat{G}$: 

\begin{enumerate}[(a)]
\item The spectrum of $\mathbf{J}[p]$ is the set of irreducible tempered representations of $G$ whose restriction to $K$ contains every class in $\mathpzc{C}_{p}$. 
\item The spectrum of $\mathbf{C}[p]$ is the set of irreducible tempered representations of $G$ whose set of lowest $K$-types is exactly $\mathpzc{C}_{p}$.
\end{enumerate}

\noindent A consequence of (a) is that the union $\bigcup \limits_{p \in \N} \mathbf{I}[p]$ is dense in $C^\star(G)$. Thus, the family  $( \mathbf{I}[p] )_{p \in \N}$ is a  filtration of ${C}^\star_{r}(G)$ with dense image.


\subsubsection{Morita-equivalence with a commutative algebra}
\noindent Fix a class $\lambda$ in $\widehat{K}$ and suppose $\pi$ is an irreducible tempered representation of $G$, acting on a Hilbert space $\mathcal{H}$, containing $\lambda$ as a lowest $K$-type. Every vector $v$ in the $K$-isotypical subspace $\mathcal{H}^\lambda$ of $\mathcal{H}$  determines a matrix element $\mathbf{p}^v_{\lambda}$, as in \eqref{choix_v}; the results of \S \ref{subsec:notations_coeffs} show that $\mathbf{p}^v_{\lambda}$  defines a \emph{projection} in the multiplier algebra of $\mathbf{C}[p]$. Because the $K$-type $\lambda$ occurs with multiplicity one in $\pi$ (Proposition \ref{multiplicite}), the induced projection $\pi(\mathbf{p}^v_{\lambda})$ of $\mathcal{H}$ has rank one. 

This makes the idempotent $\mathbf{p}^v_{\lambda}$ quite special, and gives precise information about the representation-theoretic structure of the subquotient $\mathbf{C}[p]$.

\begin{lemme}[Lemma 6.1 in \cite{Higson2008}] \label{lemme_cle} {Let $\mathbf{C}$ be a $C^*$-algebra and $\mathbf{p}$ be a projection in the multiplier algebra of  $\mathbf{C}$. Assume that for every irreducible representation $\pi$ of  $\mathbf{C}$, the operator $\pi(\mathbf{p})$ is a rank-one projection. Then }
{\begin{enumerate}[(i)] \item{ $\mathbf{CpC} = \mathbf{C}$;} \item{ $\mathbf{pCp}$ is a commutative $C^*$-algebra; } 
\item{the spectrum $\widehat{C}$ is a Hausdorff locally compact space;} \item{Let $\mathcal{C}_0(\widehat{\mathbf{C}})$ be the algebra of continuous functions on $\widehat{\mathbf{C}}$ that vanish at infinity. The map $a \mapsto \widehat{a}$ from  $\mathbf{pCp}$ to $\mathcal{C}_0(\widehat{\mathbf{C}})$ defined by \[ \forall \pi \in \widehat{\mathbf{C}}, \quad \pi(a) = \widehat{a}(\pi) \pi(\mathbf{p}) \] is an isomorphism of $C^*$-algebras.} \end{enumerate} }\end{lemme}

\subsubsection{Some precisions on the subquotient's spectrum} 

\noindent Part of Lemma \ref{lemme_cle} says that the spectrum of $\mathbf{C}[p]$ is a locally compact Hausdorff space. We can give a somewhat concrete picture for that space using the classification of irreducible tempered representations.

Some additional notations will be necessary. Fix a maximal abelian subalgebra $\mathfrak{a}$ in $\mathfrak{p}$. Choose a family $P_1, \dots,  P_r$ of nonconjugate cuspidal parabolic subgroups in $G$, with one member in the family for each conjugacy class of cuspidal parabolic subgroups; arrange the choice so that in the Langlands decomposition $M_i A_i N_i$ of $P_i$, the Lie algebra  $\mathfrak{a}_i$ is contained in $\mathfrak{a}$.  Anticipating a later need for further notation, we will write $M_i^{\mathfrak{p}}$ for $\exp_{G} \left(\mathfrak{m}_i \cap \mathfrak{p} \right)$ and recall that, with the obvious notations for the Lie algebras, the Iwasawa map from $K \times \left(\mathfrak{m}_i \cap \mathfrak{p} \right)\times \mathfrak{a}_i \times \mathfrak{n}_i$ to $G$ is a diffeomorphism.

 Now, for each $i$ in $\{1, \dots, r\}$, write $K_i$ for the maximal compact subgroup $K \cap M_i$ in $M_i$ and fix a linear ordering $\widehat{K_i} = \{\lambda_1, \lambda_2, \dots \}$ in such a way that if $\norme{\lambda_n }_{\widehat{K_i}}<\norme{\lambda_m }_{\widehat{K_i}}$, then $n < m$. By \emph{discrete parameter}, we will henceforth mean a couple $(i, n)$ with $i$ in $\{1, \dots,r\}$ and $n$ in $\N$.

\begin{lemme} \label{lemme_ktypes} { Fix an element $\mathpzc{C}$ of $\mathpzc{Classes}$. There exists a discrete parameter $(i_{0}, n_{{0}})$ and a subset  $\widehat{\mathfrak{a}}[p] $ of $\mathfrak{a}_{i_{0}}^\star$ such that}
\begin{itemize}
\item[$\bullet$] {$\mathbf{V}_{{M_{i_{0}}}}(\mu_{n_{{0}}})$ is a discrete series or nondegenerate limit of discrete series representation of $M_{i_{0}}$,}
\item[$\bullet$] {Every irreducible tempered representation of $G$ whose set of lowest K-types equals $\mathpzc{C}$ is equivalent to exactly one of the}
\[ \ind_{M_{i_{0}}A_{i_0}N_{i_0}}^{G} \left[ \mathbf{V}_{{M_{i_{0}}}}(\mu_{n_{{0}}}) \otimes e^{i\chi} \right], \chi \in \widehat{\mathfrak{a}}[p]. \]
\item[$\bullet$] {For every $\chi$ in $\widehat{\mathfrak{a}}[p]$, the parabolic subgroup $P_\chi$ of \S \ref{subsec:rappels_mackey} satisfies $P_\chi \supset P_{i_0}$ (in fact $M_\chi \supset M_{i_0}$, while $A_\chi \subset A_{i_0} $ and $N_\chi \subset N_{i_0}$).}
\end{itemize}\end{lemme}

 The second point identifies the spectrum of $\mathbf{C}[p]$ with $\widehat{\mathfrak{a}}[p]$ as a set; we shall see that when $\widehat{\mathfrak{a}}[p]$ is equipped with the topology that it inherits from Euclidean space, the identification becomes a homeomorphism.\\

\begin{exemple} When $G = SL(2, \R)$, the spherical principal series representations with nonzero continuous parameter have the same lowest $K$-type as the (irreducible) spherical principal series representation with continuous parameter zero (the $K$-type is the trivial one). The nonspherical principal series representation with nonzero continuous parameter have two distinct lowest $K$-types, and the nonspherical principal series representation with continuous parameter zero is reducible and splits into the two limits of discrete series representations. The other irreducible tempered representations are in the discrete series. So in each case, $\widehat{\mathfrak{a}}[p]$ is either a closed half-line, an open half-line or a single point.   \end{exemple}

In general $\widehat{\mathfrak{a}}[p]$ consists of an open Weyl chamber in a subspace of $\mathfrak{a}$, together with part of one of its walls, and that part-of-wall is itself stratified analogously until one reaches a minimal dimension (which might be nonzero).  \\

\noindent \emph{Proof of Lemma \ref{lemme_ktypes}.}  For every $\pi$ in $\widehat{G}$, we use double induction from the Mackey-Higson parametrization of $\widetilde{G}$ in Theorem \ref{correspondance}, the Knapp-Zuckerman classification of representations, and basic structure theory, to obtain the existence of an element $\chi$ of $\mathfrak{a}^\star$ and a discrete parameter $(i(\chi), n(\chi))$ such that 
\begin{itemize}
\item[$\bullet$] $\mathbf{V}_{{M_{i(\chi)}}}(\mu_{n(\chi)})$ is a discrete series or nondegenerate limit of discrete series representation,
\item[$\bullet$]  the representation $\pi$ is equivalent with  $\ind_{M_{i(\chi)}A_{i(\chi)}N_{i(\chi)}}^{G} \left[ \mathbf{V}_{{M_{i(\chi)}}}(\mu_{n(\chi)}) \otimes e^{i\chi} \right]$, 
\item[$\bullet$] the parabolic subgroup $P_\chi$ contains $P_{i(\chi)}$.
\end{itemize} 
Now if $\pi = \mathbf{M}(\chi, \mu)$ and $\pi'=\mathbf{M}(\chi', \mu')$ are representations in $\widetilde{G}$ which have the same lowest $K$-types, Vogan's work shows that the discrete parameters  $(i(\chi), n(\chi))$ and  $(i(\chi'), n(\chi'))$ must be the same (see e.g. \cite{AAMackey}, Lemma 4.2).  Write $(i_0, n_0)$ for the discrete parameter shared by all representations in $\widetilde{G}$ whose set of lowest $K$-types is $\mathpzc{C}$. Then every representation in $\widetilde{G}$ whose set of lowest $K$-types is $\mathpzc{C}$ is equivalent with  one of the $\ind_{M_{i_{0}}A_{i_0}N_{i_0}}^{G} \left[ \mathbf{V}_{{M_{i_{0}}}}(\mu_{n_{{0}}}) \otimes e^{i\chi} \right], \chi \in \mathfrak{a}_{i_0}^\star$.  Lemma \ref{lemme_ktypes} follows from this and from the uniqueness properties in the Knapp-Zuckerman classification. \qed

\subsubsection{An explicit formula for the isomorphism in Lemma \ref{lemme_cle}.}\label{subsec:explicite}

\noindent In this subsection, we fix a class $\mathpzc{C} = \mathpzc{C}_p$ in $\mathpzc{Classes}$ and write $(i_{0}, n_{0})$ for the discrete parameter attached to $\mathpzc{C}$ by Lemma \ref{lemme_ktypes}. We fix an element $\chi$ in $\widehat{\mathfrak{a}}[p]$; this determines a representation $\pi$ in the spectrum of $\mathbf{C}[p]$. In order to exploit the structural information provided by Lemma \ref{lemme_cle}, we will need to give a concrete form to the transform $\widehat{a}(\pi)$ that occurs there. 

 We fix \emph{one} $K$-type $\lambda_p$ in the class $\mathpzc{C}_{p}$. Since $p$ will remain fixed until the end of section \ref{sec:cstar}, we will usually remove the subscript $p$ from $\lambda$. 
 
Picking a carrier Hilbert space $\mathcal{V}^{n_0}$ for a representation $\sigma_{n_0}$ of $M_{i_0}$ whose equivalence class is  $\mathbf{V}_{{M_{i_{0}}}}(\mu_{n_{{0}}})$, write $\mathcal{H}^{\chi}_{i_{0}, n_{0}}$ for the Hilbert space carrying $\ind_{M_{i_{0}}}^{G} \left[ \mathbf{V}_{{M_{i_{0}}}}(\mu_{n_{{0}}}) \otimes e^{i\chi} \right]$ in the usual induced picture: the~completion~of 
\[\left\{ \xi: G\  \overset{\text{smooth}}{\underset{\text{comp. supp.}}{\longrightarrow}} \  \mathcal{V}^{n_0} \ | \ \xi(gman) = a^{-i\chi-\rho} \sigma_{n_0}(m)^{-1} \xi(g) \ \text{ for } \  \left( g, m,a,n \right) \in G \times M_{i_0} \times A_{i_0} \times N_{i_0} \right\} \] 
in the norm associated to the inner product $\langle \xi_1, \xi_2 \rangle = \displaystyle \int_K \langle \xi_1(k), \xi_2(k) \rangle_{\mathcal{V}^{n_0}}dk$. In the above equation, $\rho$ is the half-sum of positive roots for $(\mathfrak{g}, \mathfrak{a}_{i_0})$ in the ordering defined by $N_{i_0}$. 

Write $\pi_{n_{0} \chi}$  for the usual morphism from $G$ to the unitary group of $\mathcal{H}^{\chi}_{i_{0}, n_{0}}$  : for every $g$ in $G$ and $\xi$ in $\mathcal{H}^{\chi}_{i_{0}, n_{0}}$, $\pi_{n_{0} \chi}(g) \xi$ is the map $x \mapsto \xi(g^{-1}x)$.

Inside the carrier space $\mathcal{V}^{n_0}$ for $\sigma_{n_0}$, consider the $K_{i_0}$-isotypical subspace $W$ which corresponds to the $K_{i_0}$-type $\mu_{n_0}$. Upon decomposing the $\lambda$-isotypical $K$-invariant component in $\mathcal{H}^{\chi}_{i_{0}, n_{0}}$, say $V$, into $K_{i_{0}}$-invariant parts, Frobenius reciprocity says the $K_{i_{0}}$-type $\mu_{n_{{0}}}$ appears exactly once. Fix a $K$-equivariant identification between the corresponding $K_{i_{0}}$-irreducible subspace and $W$, write $\tilde{v}$ for the vector in the $\lambda$-isotypical subspace $V$ which the identification assigns to any $v$ in $W$. Now, choose an arbitrary nonzero $v_{\text{spec}}$ in $W$, and introduce as in \eqref{choix_v} a matrix coefficient of our irreducible $K$-module of type $\lambda$ as
\begin{equation} \label{matel} \mathbf{p}_{\lambda} = \mathbf{p}^{v_{\text{spec}}}_{\lambda} = k \mapsto \langle  \tilde{v}_{\text{spec}},\pi_{n_0 \chi}(k) \tilde{v}_{\text{spec}} \rangle \end{equation}
This is our choice of $v$ in \eqref{choix_v};  henceforth we shall assume that it is this matrix element that is chosen in the definition of $\mathbf{J}[p]$ and $\mathbf{C}[p]$. Note that the smooth function on $K$ which we just defined does not depend~on~$\chi$.


 The projection $\pi_{n_{0} \chi}(\mathbf{p}_{\lambda})$ has rank one. To make Lemma \ref{lemme_cle} explicit, we now \emph{imitate}  Higson's analysis of the complex semisimple case, and \emph{use} the calculations he made in \cite{Higson2008} for the Cartan motion group rather than complex semisimple groups.
 
  A first step is to identify the range of $\pi_{n_{0} \chi}(\mathbf{p}_{\lambda})$.  Choose a basis $\{v_\alpha \}_{i=1,\dots, d(\mu_{n_0})}$  for  $W$  in such a way that one of its vectors, say $v_{\alpha_0}$, is the $v_{\text{spec}}$ that occurs in \eqref{matel}. Define a function from $K$ to $W$ as
\begin{equation} \label{zeta} \zeta_{n_0 \chi} :  k \mapsto  \sum \limits_{\alpha = 1}^{d(\mu_{n_0})}  \langle \pi_{n_0 \chi}(k) \tilde{v}_\alpha,  \tilde{v}_{\text{spec}} \rangle v_\alpha = \frac{1}{d(\lambda)} \sum \limits_{\alpha = 1}^{d(\mu_{n_0})} e_{\alpha_0 \alpha}(k)v_\alpha, \end{equation}
where we used the notations of \S \ref{subsec:notations_coeffs} and wrote $d(\lambda)$ and $d(\mu_{n_0})$  for the dimensions of $V_\lambda$ and $W$, respectively. We can obtain a vector in the representation space $\mathcal{H}^{\chi}_{i_{0}, n_{0}}$  by extending $\zeta_{n_0 \chi}$ to $G$, setting 
\[ \xi_{n_0 \chi}(km^{\mathfrak{p}}an) = \frac{d(\lambda)^{1/2}}{\text{vol}(K)^{1/2} d(\mu_{n_0})^{1/2}} \ e^{- i\chi-\rho}(a)     \sum \limits_\alpha \langle \pi_{n_0 \chi}(k)  \tilde{v}_\alpha, \tilde{v}_{\alpha_0} \rangle \cdot \sigma_{n_0}(m^\mathfrak{p})^{-1} \left[ v_\alpha \right]. \]
\begin{lemme} {The operator $\pi(\mathbf{p}_\lambda)$ on $\mathcal{H}^{\chi}_{i_{0}, n_{0}}$ agrees with the orthogonal projection on $\xi_{n_0 \chi}$}.\end{lemme}
 This is easily proved using the formula for the action of $K$ and the inner product on the representation space, and a repeated application of the Schur-Weyl orthogonality relations. Now put
\begin{align} \label{coeff} \widehat{f}^{[p]}: \hspace{0.3cm} \widehat{\mathfrak{a}}[p] &\rightarrow \C \nonumber \\  \chi & \mapsto \int_{G} f(g) \langle \xi_{n_{0} \chi}, \pi_{n_{0}\chi}(g) \xi_{n_{0} \chi} \rangle \end{align}
as soon as $f$ is a smooth and compactly supported function on $G$. If $ \mathbf{p}_{\lambda}  \underset{K}{\star} f  \underset{K}{\star}  \mathbf{p}_{\lambda}  = f$, then $\pi_{n_{0} \chi}(f)$ is proportional to $\pi_{n_{0} \chi}(\mathbf{p}_\lambda)$, and given the definition of  $\pi_{n_{0} \chi}(f)$,  we see that $\pi_{n_{0} \chi}(f) =  \widehat{f}^{[p]}(\chi) \ \pi_{n_{0} \chi}(\mathbf{p}_{\lambda})$. This is a first step in making Lemma \ref{lemme_cle} explicit. The explicit form for $ \widehat{f}^{[p]}$ to be given below will prove that it is continuous and vanishes at infinity as a function of $\chi$, so we can summarize the above in the following statement (compare Lemma 6.10 in \cite{Higson2008}). 

\begin{prop} \label{prop_morita} \label{structure_G} {By associating, to any smooth and compactly supported function $f$ on $G$ such that $\mathbf{p}_{\lambda} \underset{K}{\star} f \underset{K}{\star}  \mathbf{p}_{\lambda} = f$, the element $ \widehat{f}^{[p]}$ of $\mathcal{C}_{0}(\widehat{\mathfrak{a}}[p])$, one obtains a $C^\star$-algebra isomorphism between $\mathbf{p}_{\lambda} \mathbf{C}[p] \mathbf{p}_{\lambda}$ and $\mathcal{C}_{0}(\widehat{\mathfrak{a}}[p])$.}\end{prop}

\noindent As promised just after Lemma \ref{lemme_ktypes}, this identifies the spectrum of $\mathbf{C}[p]$ with $\widehat{\mathfrak{a}}[p]$ homeomorphically. \\

\noindent It will be important later on to have a completely explicit formula for $ \widehat{f}^{[p]}$, so we now record a closed form for \eqref{coeff}, to be used in section \ref{defo_cstar} below. 

 Given $\alpha, \beta$ in $\{1, \dots,  d(\mu_{n_0}) \}$, define a function of  $\chi \in \widehat{\mathfrak{a}}[p]$:
\begin{equation} \label{tr} \widehat{f}^{ \ p}_{\alpha,\beta}(\chi) = \frac{\text{Vol}(K)}{d(\lambda)} \int_{N_{i_0}} dn \int_{A_{i_0}} da \ a^{i\chi+\rho} \int_{M_{i_0}^\mathfrak{p}} dm \left( f \underset{K_{i_0}}{\star} d_{\alpha_0 \alpha} \underset{K_{i_0}}{\star} d_{\beta \alpha_0} \right) (nam) \langle v_\beta,  \sigma_{n_0}(m^{-1}) \left[ v_\alpha \right] \rangle . \end{equation}
\begin{lemme} \label{transforme} {The element $ \widehat{f}^{[p]}$ of $\mathcal{C}_{0}(\widehat{\mathfrak{a}}[p])$ can be expressed as}
\[ \widehat{f}^{[p]} = \frac{1}{d(\mu_{n_0})} \sum_{\alpha, \beta = 1}^{d(\mu_{n_0})} \widehat{f}^{\ p}_{\alpha, \beta}. \]\end{lemme}
\begin{proof} We expand \eqref{coeff}, closely following the calculations on page 15 of \cite{Higson2008}. We start from the fact that $ \displaystyle \langle \xi_{n_{0} \chi}, \pi_{n_{0}\chi}(g) \xi_{n_{0} \chi} \rangle = \int_K \langle \xi_{\mu_{n_0} \chi}(k) \ , \  \xi_{\mu_{n_0} \chi}(g^{-1}k) \rangle dk $; after a change of variables $g \leftarrow g^{-1}k$, and inserting the necessary normalizations to have the $e_{\alpha \beta}$ appear, we find
\begin{align*} \widehat{f}^{[p]}(\chi) & = \int_G \left( \int_K  f(kg^{-1}) \langle \xi_{\mu_{n_0} \chi}(k), \xi_{\mu_{n_0} \chi}(g) \rangle dk \right) dg\\
 & = \int_G \left( \int_K  f(kg^{-1}) \frac{1}{\text{Vol}(K)^{1/2}d(\lambda)^{1/2} d(\mu_{n_0})^{1/2}}\sum \limits_{\alpha} \overline{e_{\alpha_0 \alpha}(k)} \langle v_\alpha, \xi_{\mu_{n_0} \chi}(g) \rangle dk \right) dg \\
 & = \left( \frac{1^{}}{\text{Vol}(K)^{}d(\lambda)^{1/2} d(\mu_{n_0})^{}} \right)^{1/2} \sum \limits_{\alpha} \int_G \text{Vol}(K) \left( e_{\alpha \alpha_0}  \underset{K}{\star} f\right)(g^{-1}) \langle v_\alpha, \xi_{\mu_{n_0} \chi}(g)\rangle dg \\ &= \frac{1}{d(\lambda)d(\mu_{n_0})} \sum \limits_{\alpha, \beta} \int_{N_\chi} dn \int_{A_\chi}da \int_{M^\mathfrak{p}_\chi}dm \int_K du \ ( e_{\alpha\alpha_0}  \underset{K}{\star} f)(n^{-1}a^{1}m^{-1} u^{-1}) e_{ \alpha_0 \beta}(u) \langle v_\alpha, \sigma(m^{-1}) v_\beta \rangle a^{-i\chi-\rho}.
 \end{align*} 
 Hence,
\begin{equation*} 
 \widehat{f}^{[p]}(\chi) =  \frac{ \text{Vol}(K)}{d(\lambda) d(\mu_{n_0})} \sum \limits_{\alpha, \beta} \int_{N_{\chi}}dn \int_{A_{\chi}} da\  a^{-i\chi-\rho} \int_{M^\mathfrak{p}_{\chi}}\left[ e_{\alpha \alpha_0}  \underset{K}{\star} f  \underset{K}{\star} e_{\alpha_0 \beta} \right](n^{-1}a^{-1}m^{-1})  \langle v_\alpha, \sigma(m^{-1}) v_\beta \rangle.\end{equation*} 
  To relate this to the quantities $\widehat{f}^{\ p}_{\alpha,\beta}$ of \eqref{tr}, we need to have the $d_{\alpha \beta}$ enter the formula in place of the $e_{\alpha \beta}$. We use the orthogonality relations \eqref{weyl} to observe that 
\[ e_{\alpha \alpha_0} \underset{K}{\star} f \underset{K}{\star}  e_{\alpha_0 \beta} = d_{\alpha \alpha_0}  \underset{K_{i_0}}{\star}  \left( e_{\alpha_0 \alpha_0} \underset{K}{\star}  f \underset{K}{\star}  e_{\alpha_0 \alpha_0} \right) \underset{K_{i_0}}{\star}  d_{\alpha_0 \beta}=d_{\alpha \alpha_0}  \underset{K_{i_0}}{\star}  \left( \mathbf{p}_{\lambda} \underset{K}{\star}  f \underset{K}{\star}  \mathbf{p}_{\lambda} \right) \underset{K_{i_0}}{\star}  d_{\alpha_0 \beta}.\] Because of our hypothesis on $f$, this is actually equal to $d_{\alpha \alpha_0}  \underset{K_{i_0}}{\star}  f \underset{K_{i_0}}{\star}  d_{\alpha_0 \beta}$.

 To get two convolutions on the right of $f$ instead of one on each side of $f$, we shorten the above formula by writing $\Gamma^{n_0}_{\alpha \beta}(m)$ for $ \langle v_\alpha, \sigma(m^{-1}) v_\beta \rangle$, and we use the structure properties of the parabolic subgroup $P_{i_0}$, along with the fact that $\chi$ is $K_\chi$-invariant (hence $K_{i_0}$-invariant because of the last assertion in Lemma \ref{lemme_ktypes}), to obtain
\begin{align*} & \int_{N_{i_0}} dn \int_{A_{i_0}} da\  a^{-i\chi-\rho} \int_{M^\mathfrak{p}_{i_0}}\left[ d_{\alpha \alpha_0} \star f \star d_{\alpha_0 \beta} \right](n^{-1}a^{-1}m^{-1}) \Gamma^{n_0}_{\alpha \beta}(m) = \\
 &  \int_{N_{i_0}} dn  \int_{M^\mathfrak{p}_{i_0}} dm \Gamma^{n_0}_{\alpha \beta}(m)\int_{K_{i_0}} dk_1 \int_{K_{i_0}} dk_2 \int_{A_{i_0}}da\  a^{-i\chi-\rho} d_{\alpha \alpha_0}(k_1) f(k_1^{-1} n^{-1}m^{-1}a^{-1}  k_2 ) d_{\alpha_0 \beta}(k_2)   = \\ 
 &   \int_{N_{i_0}} dn  \int_{M^\mathfrak{p}_{i_0}} dm \Gamma^{n_0}_{\alpha \beta}(m)\int_{K_{i_0}}  dk_1 \int_{K_{i_0}} dk_2 \int_{A_{i_0}}da\  a^{-i\chi-\rho} d_{\alpha \alpha_0}(k_1) f( n^{-1}m^{-1}a^{-1} k_1^{-1} k_2 ) d_{\alpha_0 \beta}(k_2) =\\
 &  \int_{N_{i_0}} dn  \int_{A_{i_0}} da\  a^{-i\chi-\rho} \int_{M^\mathfrak{p}_{i_0}}\left[  f \star d_{\alpha \alpha_0} \star d_{\alpha_0 \beta} \right](n^{-1}a^{-1}m^{-1}) \Gamma^{n_0}_{\alpha \beta}(m)  \end{align*}
(the stars are now convolutions over $K_{i_0}$; between the first and second line, we used the fact that $M_{i_0}$ centralizes $A_{i_0}$; between the second and third line, we used the fact that $K_{i_0}$ is contained in $M_{i_0}$, and thus leaves $A_{i_0}$ invariant and normalizes $N_{i_0}$, to perform the change of variables $m \leftarrow k_1^{-1} m k_1$, $n \leftarrow k_1^{-1} n k_1$, $a \leftarrow k_1^{-1} a k_1$). The last quantity is that which appears in \eqref{tr}: this proves the lemma. \end{proof}
\subsection{Subquotients of the motion group's algebra}\label{subsec:sous_quotient_G0}

\noindent Taking up the notations of sections \ref{subsec:notations_coeffs} and \ref{subsec:sous_quotient_G}, define a closed ideal in the reduced $C^\star$ algebra ${C}_{r}^\star(G_{0})$ by setting, as before,
\[ \mathbf{J}^0[p]  = \bigcap \limits_{\lambda \in \mathpzc{C}_{p}}  {C}_{r}^\star(G_{0})  \mathbf{p}_{\lambda} {C}_{r}^\star(G_{0}),  \]
and a subquotient of ${C}_{r}^\star(G_{0})$ by setting
\[ \mathbf{C}^0[{p}] = \left(\mathbf{J}^0[1] + \dots + \mathbf{J}^0[p]\right)/\left(\mathbf{J}^0[1]+\dots+\mathbf{J}^0[p-1] \right).\]
 As before, the spectrum of $\mathbf{C}^0[p]$ is the set of unitary irreducible representations of $G_{0}$ whose set of lowest $K$-types is exactly $\mathpzc{C}_{p}$. Because the Mackey-Higson bijection preserves lowest $K$-types (Proposition \ref{preservation_ktypes}), we know that the spectrum of $\mathbf{C}^0[p]$ can be identified \emph{as a set} with $\widehat{\mathfrak{a}}[p]$. 

The multiplicity-one property in Proposition \ref{multiplicite}(b) means that Lemma \ref{lemme_cle} is still applicable: the subquotient $\mathbf{C}^0[p]$ is thus Morita-equivalent with the algebra of continuous functions, vanishing at infinity, on its spectrum $-$ viewed as a topological space. To make that the analogy between $\mathbf{C}^0[p]$ and $\mathbf{C}[p]$ more precise, the next step is to identify the spectrum of $\mathbf{C}^0[p]$ with $\widehat{\mathfrak{a}}[p]$ not only as a set, but also as a {topological space}, writing down an analogue of Proposition \ref{prop_morita}. The explicit calculations in Higson's work are actually sufficient for this, so instead of \emph{imitating} the results of section 5.3 of \cite{Higson2008}, we shall \emph{invoke} them.

For every $\chi$ in $\widehat{\mathfrak{a}}[p]$, Lemma \ref{lemme_ktypes} and Proposition \ref{preservation_ktypes}  show that there is exactly one element $\mu_{p, \chi}$ in $\widehat{K_{\chi}}$ for which the set of lowest K-types of the $G_0$ representation $\mathbf{M}_{0}(\chi, \mu_{p, \chi})$ (defined in \S \ref{subsec:rappels_mackey}) is exactly $\mathpzc{C}_{p}$. Let us write $\widetilde{W}$ for a carrier vector space of $\mu_{p, \chi}$.

Viewing $\chi$ as an element of $\mathfrak{p}^\star$ which vanishes on the orthogonal of $\mathfrak{a}_\chi$, we realize $\mathbf{M}_{0}(\chi, \mu_{p, \chi})$ as the completion $\mathcal{H}_{0}$ of 
\[ \left\{ \xi: G_0  \overset{\text{smooth}}{\underset{\text{comp. supp.}}{\longrightarrow}}\ \widetilde{W} \ \ \big| \ \  \xi(gkx) = \mu_{p, \chi}(k)^{-1} \chi(x)^{-1} \xi(g) \text{ for  }(k, x, g)  \in K_\chi \times \mathfrak{p} \times G_0\right\}\]
 in the norm induced by the scalar product between restrictions to $K$. We write $\pi^0_{\chi, \mu}$ for the $G_0$-action on induced by left translation.

At this point, we cast a look backwards to the situation for reductive groups. We can view the $\lambda_p$-isotypical subspace $V$ of the $G$-representation $\mathcal{H}^{\chi}_{i_0,n_0 }$ as a $K$-module, and restrict it to $K_\chi$; it must then contain the $K_\chi$-type $\mu$ with multiplicity one:  the class of $\pi_{n_{0}, \chi}$ in $\widetilde{G}$ is the representation $\mathbf{M}(\chi, \mu_{p, \chi})$ defined in \S \ref{subsec:rappels_mackey}, and the argument at the beginning of \S \ref{subsec:explicite} can be repeated. This makes it possible to pin down a distinguished unit vector $\zeta_{p, \chi}$ in $\mathcal{H}_0$ by copying the definition of $\zeta_{n_0 \chi}$ (with the $K$-module $W$ of \S \ref{subsec:explicite} replaced by the current, isomorphic, $K$-module $\widetilde{W}$).

Recall from the multiplicity-one property in Proposition \ref{multiplicite} that the projection $\pi^0_{\chi, \mu}(\mathbf{p}_{\lambda_{}})$ has rank one. Higson proved in \cite[Lemma 5.8]{Higson2008} that the range of that orthogonal projection is the line generated by $\zeta_{p, \chi}$.~If~we~set
\begin{equation} \label{fourier} \widehat{f}(\chi) = \frac{ \text{vol}(K) }{d(\lambda)} \int_{\mathfrak{p}} f(x) e^{i\langle \chi, x\rangle} dx \end{equation}
when $f$ is a smooth  function on $G_{0}$ with compact support, then the condition $\mathbf{p}_{\lambda} \underset{K}{\star} f \underset{K}{\star}  \mathbf{p}_{\lambda} = f$ leads to the equality
\[ \pi^0_{\chi, \mu}(f) =  {\widehat{f}(\chi)} \ \pi^0_{\chi, \mu}(\mathbf{p}_{\lambda_{}}); \]
see \cite{Higson2008},  Lemma 6.12. 

Now, we know from Lemma \ref{lemme_cle} that if we view $\widehat{\mathfrak{a}}[p]$ as the spectrum of $\mathbf{C}^0[p]$ and equip it with the corresponding (locally compact and Hausdorff) Fell topology, then $\widehat{f}$ becomes a continuous function of $\chi$ that vanishes at infinity. But of course $\widehat{f}$ is  the ordinary Fourier transform of $f$ (up to a multiplicative constant), so  it is also a continuous function on $\widehat{\mathfrak{a}}[p]$ when the latter is equipped with the topology inherited from that of Euclidean space. We can thus summarize the situation with the following statement. 

\begin{prop} \label{transform_G0} {By associating, to any smooth and compactly supported function $f$ on $G_{0}$ such that $\mathbf{p}_{\lambda} \underset{K}{\star} f \underset{K}{\star} \mathbf{p}_{\lambda}$, the element $\widehat{f}$ of $\mathcal{C}_{0}(\widehat{\mathfrak{a}}[p])$, one obtains a $C^\star$-algebra isomorphism between $\mathbf{p}_{\lambda} \mathbf{C}^0[p] \mathbf{p}_{\lambda}$ and $\mathcal{C}_{0}(\widehat{\mathfrak{a}}[p])$.}\end{prop} 

We note that the transform at $\chi$ involves, on the $G_0$-side, the representation $\pi^0_{\chi, \mu}$, and on the $G$-side the representation $\pi_{n_0 \chi}$. These two do correspond to one another in the Mackey-Higson bijection, so Propositions \ref{prop_morita} and \ref{transform_G0} together have the following consequence.

\begin{cor} Fix a class $\mathpzc{C} \in \mathpzc{Classes}$ of lowest $K$-types, write $\widetilde{G}[\mathpzc{C}]$  $($resp. $\widehat{G_0}[\mathpzc{C}])$ for the set of irreducible tempered representations of $G$ (resp. unitary irreducible representations of $G_0$) whose set of lowest $K$-types is exactly $\mathpzc{C}$. The Mackey-Higson bijection of Theorem \ref{correspondance} induces a homeomorphism between $\widehat{G_0}[\mathpzc{C}]$ and $\widetilde{G}[\mathpzc{C}]$. \end{cor}

\section{Deformation of the reduced $C^\star$-algebra and subquotients} \label{defo_cstar}

We now come back to the contraction family $(G_t)_{t \in \R}$. For every $t>0$, going through the constructions of \S \ref{sec:cstar} yields a $G_t$-invariant analogue of Proposition 3.4. Does this furnish a way to view Propositions \ref{transform_G0} and \ref{prop_morita} as deformations of one another?

There is an obvious difference between the transforms that occur in the two results. On the reductive side, the explicit formula in Lemma \ref{transforme} involves, in a somewhat delicate way, the behavior away from $K$ of some matrix elements of the discrete series (or limit of discrete series) representation $\sigma_{n_0}$  from \S \ref{subsec:explicite}. On the Cartan motion group side, a striking (and perhaps unexpected) feature of \S \ref{subsec:sous_quotient_G0} is that the transform in Proposition \ref{transform_G0} does not involve any contribution of the discrete parameter $n_0$ (except for the global dimension factor $d(\lambda)$).  
 
In spite of this difference, we will see that the transforms do deform onto one another as the contraction is performed. The rigidity phenomenon in Theorem A$(ii)$ thus appears to be related to the collapse of the ``discrete-series-away-from-$K$'' terms as we perform the contraction to the Cartan motion group. Although we will need only elementary arguments below, this collapse does call to mind the phenomena occurring at the level of geometric realizations that we studied in \cite{AAContractions} (see \S 3.1 there for the discrete series case).  As announced, we will here again follow \cite{Higson2008} closely (see \S 6.2 and 6.3 there).

\subsection{Subquotients of the continuous field and their spectra}

\noindent  Consider the continuous field $\mathcal{C}$ from \eqref{champ_c}; as before, for all $p$ in $\N$, we let $\mathcal{J}[p]$ be the ideal $\bigcap \limits_{\lambda \in \mathpzc{C}_{p}} \mathcal{C} \mathbf{p}_{\lambda} \mathcal{C}$ of $\mathcal{C}$, set $\mathcal{I}[p]=\mathcal{J}[1] + \dots \mathcal{J}[p]$, so that the family $(\mathcal{I}[p])_{p \in \N}$ is a dense filtration of $\mathcal{C}$. We fix some $p$ in $\N$ (and accordingly, a discrete parameter $(i_0, n_0)$) and study the subquotient $\mathcal{C}[p] =\mathcal{I}[p+1]  \big/ \mathcal{I}[p]$ of $\mathcal{C}$.

The spectrum of  $\mathcal{C}[p]  $ can be identified with $\widehat{\mathfrak{a}}[p] \times [0,1]$ as a set. Indeed, the algebra $\mathpzc{Z}$ of continuous functions on the closed interval $[0,1]$ lies in the center of the multiplier algebra of $\mathcal{C}[p]$, so that if $\mathpzc{Z}(t)$ is the subalgebra of functions which vanish at $t$, the spectrum of $\mathcal{C}[p]$ can be identified as a set with the disjoint union over $t$ of the spectra of the quotient algebras $\mathcal{C}[p] /(\mathpzc{Z}(t)\mathcal{C}[p] )$ $-$ with each $t$, we need only associate the subset of $\widehat{\mathcal{C}[p]}$ whose elements are the representations that restrict to zero on $\mathpzc{Z}(t)\mathcal{C}[p] $ .  But for fixed $t$, the algebra $\mathcal{C}[p] /(\mathpzc{Z}(t)\mathcal{C}[p] )$ is isomorphic with the subquotient of the group algebra $C^\star_{r}(G_t)$ that corresponds to the class $\mathpzc{C}_p$, and the spectrum of that subquotient can be identified with $ \widehat{\mathfrak{a}}[p] $; this proves our claim.  

 To proceed further, we need to see how  the spectrum of  $\mathcal{C}[p]  $ can be identified with $\widehat{\mathfrak{a}}[p] \times [0,1]$ as a topological space. Suppose $f$ is a smooth and compactly supported function on the manifold $\mathpzc{G}$ of \S \ref{subsec:defo}. For each $t>0$, collect the the ingredients for a $G_t$-equivariant version of the transform that appears in Lemma \ref{transforme}:
  
 \begin{enumerate}[(i)]
\item Write $f_{t}$ for the restriction of $f$ to a smooth and compactly supported function on $G_{t}$.  Write $M_{i_{0},t} A_{i_{0},t} N_{i_{0}, t}$ for the cuspidal parabolic subgroup subgroup $\varphi_{t}^{-1} P_{i_0}$ of $G_{t}$, and note that it comes with an ordering on the $\mathfrak{a}_{i_{0},t}$-roots and an associated half-sum $\rho_t$ of positive roots.  In fact, we have $\rho_t=t\rho$ (see \cite{AAContractions}, Lemma~3.12).
 
\item Suppose $\sigma^t_{n_0}$ is an irreducible representation of $M_{i_{0}, t}$ with equivalence class $\mathbf{V}_{M_{i_{0}, t}}(\mu_{n_{0}})$; fix a carrier space $\mathcal{V}_t^{n_0}$ for $\sigma^t_{n_0}$.  We choose a basis $(v_{a}^t)_{a\in \{1, \dots,\ d(\mu_{n_0})\}}$ for the $K_{i_0}$-isotypical subspace $W_t$ of $\mathcal{V}_t^{n_0}$ that carries the lowest $K_{i_0}$-type of $\sigma^t_{n_0}$, as follows. The composition $\sigma^t_{n_0} \circ\varphi_t^{-1}$ defines an irreducible tempered representation of $M_{i_0}$ which has real infinitesimal character and lowest $K_{i_0}$-type $\mu_{n_0}$; therefore it is equivalent with $\sigma_{n_0}$, and there exists a unitary operator $\mathbf{c}_t: \mathcal{V}^{n_0} \to \mathcal{V}_t^{n_0}$ that intertwines $\sigma_{n_0}$ and  $\sigma^t_{n_0} \circ\varphi_t^{-1}$. Our choice of basis $(v_{a}^t)$ for $W_t$ is the image under $\mathbf{c}_t$ of the basis $(v_a)_{a\in \{1, \dots,\ d(\mu_{n_0})\}}$  for $W \subset \mathcal{V}^{n_0}$ which we used in \eqref{tr}.
\end{enumerate}
 
%
%

\noindent We can now gather the transforms considered in \S   \ref{subsec:sous_quotient_G} and \S \ref{subsec:sous_quotient_G0}, and define, for every $\chi$ in $\widehat{\mathfrak{a}}[p]$,
 \begin{equation} \label{transf} \hspace{-0cm} \widehat{f}^{\ p}_{a, b}(\chi, t) = \begin{cases} \displaystyle \ \  \delta_{a, b} \  \int_{\mathfrak{p}} f_{0}(x) \chi(x) dx \ \quad \quad \text{($\delta_{a,b}=$ Kronecker delta)} & \text{ if } t = 0, \\  \\ \displaystyle   \int_{  N_{i_0, t} } dn_{t} \int_{ A_{i_0,t}} da_{t} \ a_{t}^{ i \chi + \rho_{t}} \int_{ M^\mathfrak{p}_{i_0,t}} dm_{t} \left( f_{t} \star d_{\alpha_0 a} \star d_{b \alpha_0} \right) \left( n_{t}a_{t}m_{t} \right) \langle v^t_b,  \sigma^t_{n_0}(m_{t}^{-1}) \left[ v_a^t \right] \rangle & \text{if } t \neq 0.   \end{cases} \end{equation}
Last, we define 
\begin{equation} \label{sm} \widehat{f}^{p}: (\chi, t) \mapsto \sum \limits_{a, b = 1}^{d(\mu_{n_{0}})} \widehat{f}^{\ p}_{a, b}(\chi, t),\end{equation}
a map from $\widehat{\mathfrak{a}}[p] \times [0,1]$ to $\C$.\\

The rigidity property of the subquotient $\mathcal{C}[p]$ announced in Theorem A$(ii)$ boils down to the the fact that this transform has no singularity at $t=0$; that is our next result. 

\begin{lemme}\label{isom_famille} {The transform $\widehat{f}^p$  of $f$ is a continuous function on $\widehat{\mathfrak{a}}[p] \times [0,1]$.}\end{lemme}

\begin{proof} We will apply the usual theorems on the continuity of parameter integrals, but must deal with the fact that the subgroups over which the integrals are taken depend on $t$. In order to write down \eqref{transf} as an integral over a space which does not depend on $t$, we set $\mathfrak{y}_{i_{0}}$ for $\mathfrak{p} \cap \left( \mathfrak{a}_{i_{0}} \oplus  \mathfrak{m}^\mathfrak{p}_{i_{0}}\right)^\perp$ (the orthogonality is with respect to a bilinear form like that in \S \ref{subsec:rappels_mackey}); we call in the Cartan involution $\theta$ of $\g$ with fixed-point-set $\mathfrak{k}$, remark that $\theta$  is also the Cartan involution of $\mathfrak{g}_{t}$ with fixed-point-set $\mathfrak{k}$;  we then write, for $t>0$,
\[ \beta_{t}: \mathfrak{y}_{i_{0}} \rightarrow \mathfrak{n}_{i_{0}, t} \]
for the inverse of $n \mapsto n - \theta(n)$. Using exponential coordinates for $\exp_{G_{t}}(\mathfrak{p})$, remarking that $\mathfrak{m}_{i_{0}}^\mathfrak{p}$ and $\mathfrak{a}_{i_0}$ are $d\varphi_{t}(0)$-invariant, and inserting the fact that $\rho_t$ is none other than $t\rho$, we can rewrite $ \widehat{f}^{\ p}_{a, b}(\chi, t)$ as
\[\small \int_{  \mathfrak{y}_{i_0} } dy \int_{ \mathfrak{a}_{i_0}} da \ e^{i \langle \chi,a\rangle} \ e^{\langle\rho, ta\rangle} \int_{ \mathfrak{m}^\mathfrak{p}_{i_0}} dm \left( f \star d_{\alpha_0 a} \star d_{b \alpha_0} \right) \left(\exp_{{G_t}}({\beta_t y}) \exp_{G_t}(a) \exp_{G_t}(m) \right) \langle v^t_b,  \sigma^t_{n_0}(\exp_{G_t}(m)^{-1}) \left[ v^t_a \right] \rangle \normalsize\]
for $t \neq 0$.
We now remark that $f$ is a smooth, compactly supported function on $\mathpzc{G}$ if and only if there is a smooth and compactly supported function $F$ on $K \times \mathfrak{m}^\mathfrak{p}_{\chi} \times \mathfrak{a}_{\chi} \times \mathfrak{y}_{\chi} \times \R$ such that
\[ F(k, m, a, y) = 
\begin{cases}
f_{0}(k, m+a+y) \text{ if } t = 0, \\
f_{t}\left(k \exp_{G_{t}}(m) \exp_{G_{t}}(a) \exp_{G_{t}}(\beta_{t} y) \right) \text{ if } t \neq 0,
\end{cases}
\]
where $f_{t}$ is the restriction of $f$ to $G_{t}$ (compose the diffeomorphism \eqref{variables} with the proof of Lemma 6.17 in \cite{Higson2008}).
 Once we insert this, as well as the relationship between $\sigma^t_{n_{0}}$ and $\sigma_{n_0}$ detailed in point (ii) above, into \eqref{transf}, we can rewrite $\widehat{f}^{\ p}_{a, b}(\chi, t)$, for every $t>0$, as
\[ \hspace{-1cm} \int_{  \mathfrak{y}_{i_0} } dy \int_{ \mathfrak{a}_{i_0}} da  e^{i \langle \chi,a\rangle} \ e^{\langle\rho, ta\rangle} \int_{ \mathfrak{m}^\mathfrak{p}_{i_0}} dm \left[ F \star d_{\alpha_0 a} \star d_{b \alpha_0} \right](1, m, a, y, t)  \langle v_b,  \sigma_{n_0}(\exp_{G}(-t m)) \left[ v_a \right] \rangle. \]
To check that \eqref{transf} does define a continuous function as soon as $f$ induces an element of $\mathbf{C}[p]$, it is now enough to remember that $\mathfrak{p} = \mathfrak{m}^{\mathfrak{p}}_{i_{0}} + \mathfrak{a}_{i_{0}} + \mathfrak{y}_{i_{0}}$ and to apply the Lebesgue (dominated convergence) theorem, using the fact that $(v_{a})_{a=1\dots\dim(W)}$ is an orthonormal basis for $W$, and the Schur-Weyl relations on the matrix elements $d_{\alpha\beta}$. \end{proof}
 
Combining Lemma \ref{isom_famille} and the results in \S \ref{sec:sous_quotients}, we obtain the expected rigidity statement:

\begin{prop} \label{structure_champ} {If $f$ is a smooth and compactly supported function on $G_{0}$ and if $\mathbf{p}_{\lambda} f \mathbf{p}_{\lambda} = f$, then its transfom $\widehat{f}^{\ p}$ is a continuous function on $\widehat{\mathfrak{a}}[p] \times [0,1]$, and it vanishes at infinity. This determines an isomorphism of $C^\star$-algebras}
\[ \mathbf{p}_{\lambda} \mathcal{C}[p] \mathbf{p}_{\lambda} \rightarrow \mathcal{C}_{0}\left(\widehat{\mathfrak{a}}[p] \times [0,1] \right).\]
\end{prop}

\subsection{Proof that the Connes-Kasparov map \eqref{ConnesKasparov} is an isomorphism}

\noindent We can now prove that the map $\alpha_1$ in Theorem \ref{strategie_CK} is an isomorphism. The argument is not only close to, it is identical with, that in \S 7 of \cite{Higson2008}, so  the contents of this section are due to~Nigel~Higson. %

 Because the increasing filtration $(\mathcal{I}[p])_{p \in \N}$ has dense image in $\mathcal{C}$ and because $K$-theory commutes with direct limits, we need only prove that evaluation at $t = 1$ yields, for every $p$, an isomorphism between $\mathpzc{K}(\mathcal{I}[p])$ and $\mathpzc{K}(\mathbf{I}[p])$. By standard cohomological arguments, that conclusion will be attained if we prove that evaluation at $t=1$ induces for all $p$ an isomorphism between $\mathpzc{K}(\mathcal{C}[p])$ and $\mathpzc{K}(\mathbf{C}[p])$. 

 At this point, we recall that for fixed $p$ and $\lambda \in \mathpzc{C}_p$,  the algebras $\mathbf{p}_{\lambda}\mathbf{C}[p]\mathbf{p}_{\lambda}$ and $\mathbf{C}[p]$ are Morita-equivalent; therefore the inclusion from $\mathbf{p}_{\lambda}\mathbf{C}[p]\mathbf{p}_{\lambda}$ to $\mathbf{C}[p]$ induces an isomorphism in K-theory. The problem then reduces to showing that evaluation at $t=1$ induces an isomorphism between $\mathpzc{K}(\mathbf{p}_{\lambda}\mathcal{C}[p]\mathbf{p}_{\lambda})$ and $\mathpzc{K}(\mathbf{p}_{\lambda}\mathbf{C}[p]\mathbf{p}_{\lambda})$.

We then insert Propositions \ref{structure_champ} and  \ref{structure_G} to find that it is enough to prove that evaluation at $t=1$ induces an isomorphism between $\mathpzc{K}(\widehat{\mathfrak{a}}[p]\times[0,1])$ and $\mathpzc{K}(\widehat{\mathfrak{a}}[p])$; but that follows from the homotopy invariance of K-theory. With this the proof that $\alpha_1$ is an isomorphism is complete; the Connes-Kasparov ``conjecture'' follows.


\bibliography{mackey}
\bibliographystyle{afgou_bib}

\end{document}